\documentclass[preprint,review]{elsarticle}

\usepackage{lineno,hyperref}
\usepackage{multirow}
\modulolinenumbers[5]
\usepackage{amssymb,amstext,amsmath}
\usepackage{todonotes}
\usepackage{graphicx}
\usepackage{xcolor}

\journal{Mathematics and Computers in Simulation}

\begin{document}

\begin{frontmatter}

\title{A computational approach to extreme values and related hitting probabilities in level-dependent quasi-birth-death processes
}

\author[label1]{A. Di Crescenzo}
\ead{adicrescenzo@unisa.it}
\ead[url]{https://docenti.unisa.it/005305/home}

\author[label2]{A. G\'{o}mez-Corral\corref{cor1}}
\ead{agcorral@ucm.es}
\ead[url]{blogs.mat.ucm.es/agomez-corral}
\cortext[cor1]{Corresponding author}

\author[label3]{D. Taipe}
\ead{dtaipe@ucm.es}

\address[label1]{Dipartimento di Matematica, Universit\`{a} degli Studi di Salerno, Via Giovanni Paolo II, 132, Fisciano, I-84084, Italy}
\address[label2]{Departamento de Estad\'{\i}stica e Investigaci\'{o}n Operativa, Facultad de Ciencias Matem\'{a}ticas, Universidad Complutense de Madrid, Plaza de Ciencias 3, 28040 Madrid, Spain}
\address[label3]{Departamento de Estad\'{\i}stica y Ciencia de Datos, Facultad de Estudios Estad\'{\i}sticos, Universidad Complutense de Madrid, Avda. Puerta de Hierro s/n, 28040 Madrid, Spain}





\begin{abstract}
This paper analyzes the dynamics of a level-dependent quasi-birth-death process ${\cal X}=\{(I(t),J(t)): t\geq 0\}$, i.e., a bi-variate Markov chain defined on the countable state space $\cup_{i=0}^{\infty} l(i)$ with $l(i)=\{(i,j) : j\in\{0,...,M_i\}\}$, for integers $M_i\in\mathbb{N}_0$ and $i\in\mathbb{N}_0$, which has the special property that its $q$-matrix has a block-tridiagonal form. Under the assumption that the first passage to the subset $l(0)$ occurs in a finite time with certainty, we characterize the probability law of $(\tau_{\max},I_{\max},J(\tau_{\max}))$, where $I_{\max}$ is the running maximum level attained by process ${\cal X}$ before its first visit to states in $l(0)$, $\tau_{\max}$ is the first time that the level process $\{I(t): t\geq 0\}$ reaches the running maximum $I_{\max}$, and $J(\tau_{\max})$ is the phase at time $\tau_{\max}$. Our methods rely on the use of restricted Laplace-Stieltjes transforms of $\tau_{\max}$ on the set of sample paths $\{I_{\max}=i,J(\tau_{\max})=j\}$, and related processes under taboo of certain subsets of states. The utility of the resulting computational algorithms is demonstrated in two epidemic models: the SIS model for horizontally and vertically transmitted diseases; and the SIR model with constant population size.
\end{abstract}


\end{frontmatter}


\section{Introduction}
\label{Sect:1}
In this work, we are interested in continuous-time Markov chains ${\cal X}=\{(I(t),J(t)): t\geq 0\}$ defined on a countable state space ${\cal S}$, which we partition as $\cup_{i=0}^{\infty} l(i)$ with subsets $l(i)=\{(i,j) : j\in\{0,...,M_i\}\}$, for integers $M_i\in\mathbb{N}_0$ and $i\in\mathbb{N}_0$. The $q$-matrix of ${\cal X}$ is assumed to be conservative and irreducible, and exhibits the structured form
\begin{eqnarray}
  Q &=& \left(\begin{array}{cccccc}
      Q_{0,0} & Q_{0,1} &         & & & \\
      Q_{1,0} & Q_{1,1} & Q_{1,2} & & & \\
              & \ddots  & \ddots  & \ddots & & \\
              &         & Q_{i,i-1} & Q_{i,i} & Q_{i,i+1} & \\
              &         &           & \ddots  & \ddots    & \ddots
  \end{array}\right),
  \label{eq:1}
\end{eqnarray}
where the sub-matrices $Q_{i,i'}$ have entries $q_{(i,j),(i',j')}$, for pairs $(i,j)$ with $i\in\mathbb{N}_0$ and $j\in\{0,...,M_i\}$, and $(i',j')$ with $i'\in\{\max\{0,i-1\},i,i+1\}$ and $j'\in\{0,...,M_{i'}\}$. This means that process ${\cal X}$ can move in one step from state $(i,j)$, for $i\in\mathbb{N}_0$ and $j\in\{0,...,M_i\}$, to state $(i',j')$ only if $|i'-i|\leq 1$ and $j'\in\{0,..., M_{i'}\}$, and, consequently, the values of $I(t)$ can only change by $-1$, $0$ and $1$ units at jump times of ${\cal X}$. Process ${\cal X}$ is commonly referred to as a level-dependent quasi-birth-death (LD-QBD) process, with the first variable $I(t)$ termed the \emph{level}, and the second one $J(t)$ the \emph{phase}; for convenience, we also refer to the subset $l(i)$ as the $i$th \emph{level}. Markov chains with LD-QBD structure can be thought of as the natural generalization of the uni-variate birth-death process (see, e.g., Anderson \cite[Chapter 8]{Anderson1991}), competition processes (Iglehart \cite{Iglehart1964}; Reuter \cite{Reuter1961}) and related models, including prey-predator processes (Hitchcock \cite{Hitchcock1986}; Ridle-Rowe \cite{RidlerRowe1988}), two-species competition processes (G\'{o}mez-Corral and L\'{o}pez Garc\'{\i}a \cite{GC2012}; Ridler-Rowe \cite{RidlerRowe1978}), and birth/birth-death processes with biological applications (Ho et al. \cite{Ho2018}).
\par LD-QBD processes have been extensively investigated by several researchers, among whom are Baumann and Sandmann \cite{Baumann2010,Baumann2013}, Bright and Taylor \cite{Bright1995}, Phung-Duc et al. \cite{Phung-Duc2010}, and Takine \cite{Takine2022}. Readers are also referred to Chapter 12 of Latouche and Ramaswami \cite{Latouche1999}, the overviews of Bright and Taylor \cite{Bright1997}, and Kharoufeh \cite{Kharoufeh2011}, and references therein. The focus in these references is mainly on the stationary vector of a regular process ${\cal X}$ and its expression in terms of \emph{rate} matrices $\{R_i: i\in\mathbb{N}_0\}$, where $R_i$ consists of expected rates of visit to state $(i+1,j')$ in level $l(i+1)$ per unit of the local time of $l(i)$ spent in $(i,j)$. To be concrete, the solution of Bright and Taylor \cite{Bright1995} is developed as an extension of previous matrix-geometric methods of Neuts \cite{Neuts1981} and the logarithmic reduction solution in Ref. \cite{Latouche1993} ---which are both related to a level-independent QBD process--- for the level-dependent case. From the probabilistic interpretation of $R_i$ and a suitable continued fraction representation, Phung-Duc et al. \cite{Phung-Duc2010} provide a novel algorithmic procedure to evaluate the underlying rate matrices, which is found to be more efficient than that in Ref. \cite{Bright1995}; for a related work, see Baumann and Sandmann \cite{Baumann2010}. In Ref. \cite{Baumann2013}, Baumann and Sandmann present a matrix-analytic solution for numerically computing stationary expectations for long-run averages of additive functionals, such as long-run average costs and reward rates, without at first computing the stationary vector of ${\cal X}$. The approach of Baumann \cite{Baumann2020} (see also Baumann and Sandmann \cite{Baumann2013}) relies on truncating the state space ${\cal S}$ at a suitably selected level $l(N)$, which results in inevitable errors for which lower and upper bounds on stationary expectations are efficiently computed; see Ref. \cite[Section 4]{Baumann2020}. QBD processes on a finite state space have been extensively studied in the literature, specifically in Refs. \cite{Akar2000,DeNitto1996,Gaver1984,GC2018,GC2020}, among others. In the setting of LD-QBD processes with an \emph{explosive} state space (i.e., as the cardinality $1+M_i$ of level $l(i)$ increases exponentially with $i$), Takine \cite{Takine2022} presents approximation techniques for computing performance measures. For the analysis of quasi-stationary regime, see the work of Bean et al. \cite{Bean2000}, and Ramaswami and Taylor \cite{Ramaswami1996}.
\par In this paper, the focus is on the extreme values of the LD-QBD process ${\cal X}$ before the first entrance to states in $l(0)$, and hitting probabilities of the underlying phase process $\{J(t): t\geq 0\}$ as the maximum level of ${\cal X}$ is visited for the first time. A first motivation to study extreme values and related hitting probabilities in a LD-QBD process comes from applications. For instance, the structured form of the $q$-matrix in (\ref{eq:1}) is often obtained directly in the case of retrial queues (Artalejo and G\'{o}mez-Corral \cite{Artalejo2008}), in such a way that the length of a busy period is equivalent to the time to reach level $l(0)$ in a suitably defined LD-QBD process; for a variety of retrial queues, see the papers by Dayar and Orhan \cite{Dayar2016}, Jeganathan et al. \cite{Jeganathan2019}, Liu and Zhao \cite{Liu2010}, Phung-Duc and Kawanishi \cite{PhungDuc2014}, Phung-Duc \cite{PhungDuc2015}, Sanga and Charan \cite{Sanga2023}, and Saravanan et al. \cite{Saravanan2023}, among others. We shall demonstrate here that our approach (Section \ref{Sect:2}) does not need to exploit further properties of the retrial queue under consideration, except for the specific expressions of sub-matrices $Q_{i,i'}$, for $i\in\mathbb{N}_0$ and $i'\in\{\max\{0,i-1\},i,i+1\}$, in Eq. (\ref{eq:1}) defining the model. This is also applicable to queueing models with impatient customers and flexible matching mechanism (Chai et al. \cite{Chai2022}). In epidemics (Baumann and Sandmann \cite{Baumann2016}; Lef\`{e}vre and Simon \cite{Lefevre2020}), the first visit of ${\cal X}$ to level $l(0)$ usually amounts to the duration of an outbreak, and the maximum level of ${\cal X}$ and the state of $\{J(t): t\geq 0\}$ at the first access to $l(0)$ are related to the maximum number of simultaneously infected hosts during the outbreak and the state of the population at the peak of infection, respectively. As in the case of retrial queues, our approach (Section \ref{Sect:2}) should be seen to be a unified procedure in the analysis of epidemic models at the first time of maximum incidence of the pathogen, including multi-strain models (Chalub et al. \cite{Chalub2023}), stochastic chemical kinetics (Dayar et al. \cite{Dayar2011}), gene family evolution (Diao et al. \cite{Diao2020}), resource-limited environments (G\'{o}mez-Corral et al. \cite{GC2023c}), and recurrence in tumors (Santana et al. \cite{Santana2019}).
\par As a second motivation, let us recall that extreme values of the LD-QBD process ${\cal X}$ in (\ref{eq:1}) are analyzed by Mandjes and Taylor \cite{Mandjes2016} (see also Javier and Fralix \cite{Javier2023}) in terms of the \emph{running maximum level} $\bar{I}(T)=\sup\{I(t): t\in [0,T]\}$ attained by ${\cal X}$ at an independent exponentially distributed time $T$ and, via Erlangization, at a predetermined time $T$; for a related work in retrial queues, see G\'{o}mez-Corral and L\'{o}pez Garc\'{\i}a \cite{GC2014}. We shall complement here the work of Mandjes and Taylor \cite{Mandjes2016} by characterizing the running maximum level of a LD-QBD process at time $T$ that amounts to the first time process ${\cal X}$ makes a transition to states contained in $l(0)$, and a detailed description of the phase process $\{ J(t): t\geq 0\}$ as the running maximum is first attained. Interested readers in seeing a simpler description in epidemic models are referred to Amador et al. \cite{Amador2019}.
\par This paper proceeds as follows. In Section \ref{Sect:2}, we first introduce the basic notation for first-passage times of process ${\cal X}$, and give the formal definition of extreme value and phase as the maximum level is first attained. The distributional properties of the first time to reach the maximum level are analyzed in terms of taboo Laplace-Stieltjes transforms of certain suitably defined first-passage times (Section \ref{Subect:2.2}), and the probability law of the maximum level visited by process ${\cal X}$ (Section \ref{Subect:2.3}). In Section \ref{Sect:3}, the computational approach is exemplified with two epidemic models: the first one is used to illustrate the solution in the case of a LD-QBD process with infinitely many states (Section \ref{Subsect:3.1}), and the second one to do so in the case of a finite QBD process (Section \ref{Subsect:3.2}). Section \ref{Sect:4} concludes the paper with a brief discussion.
\section{Extreme values of a LD-QBD process and related hitting probabilities}
\label{Sect:2}
Throughout this section, the LD-QBD process ${\cal X}$ is assumed to be regular and we study the dynamics of ${\cal X}$ before the first visit to level $l(0)$, under the assumption that the first passage from any initial state $(i_0,j_0)\in {\cal S}\setminus l(0)$ to level $l(0)$ occurs in a finite time with certainty. This assumption amounts to the fact that the null column vector is the largest non-negative solution bounded above by 1 to $Q^* v=0$, where $Q^*$ is the matrix obtained by removing from $Q$ in (\ref{eq:1}) the rows and columns corresponding to states in $l(0)$; see Kulkarni \cite[Theorem 6.13]{Kulkarni2017}.
\par For ease of presentation, we adopt the notation and terminology used in our recent work \cite{GC2018,GC2020,GC2023a} and related papers \cite{Bright1995,Javier2023,Mandjes2016}. More particularly, we define the first-passage time to level $l(i)$ as $\tau_{l(i)}=\inf\{t\geq 0 : I(t)=i\}$, for integers $i\in\mathbb{N}_0$, and we let $I_{\max}$ denote the running maximum level attained by ${\cal X}$ before the first entrance into any state in level $l(0)$, and $\tau_{\max}$ be the first time to reach the value $I_{\max}$. As a result, $J(\tau_{\max})$ is the phase variable at time $\tau_{\max}$, and $I_{\max}=I(\tau_{\max})$.
\\ \\
{\bf Remark 1.}
{\it Readers should note that Javier and Fralix \cite{Javier2023} analyze the joint probability law of both the state $(I(t),J(t))$ of process ${\cal X}$ and its associated running maximum level $\bar{I}(t)$ at any arbitrary time $t$, whereas the state $(I_{\max},J(\tau_{\max}))$ is inherently linked to the random time $\tau_{\max}$, provided that the entrance of ${\cal X}$ to any state in level $l(0)$ is certain.}
\\
\par For later use, we introduce some notation as follows: $0_a$ and $1_a$ are the column vectors of order $a$ such that all their entries are equal to $0$ and $1$, respectively; $e_a(b)$ denotes the column vector of order $a$ with its $b$th entry equals $1$, and $0$ elsewhere; $0_{a\times b}$ and $I_a$ are the null matrix of dimension $a\times b$ and the identity matrix of order $a$, respectively; $^T$ denotes transposition; and $\delta_{a,b}$ denotes the Kronecker's delta.
\subsection{Statement of the problem}
\label{Subect:2.1}
Our objective is to determine the probability law of $\tau_{\max}$ on those sample paths of process ${\cal X}$ satisfying $\{I_{\max}=i,J(\tau_{\max})=j\}$, for any state $(i,j)\in\cup_{k=i_0}^{\infty}l(k)$ and a predetermined initial state $(i_0,j_0)\in {\cal S}\setminus l(0)$. To start with, we notice that
\begin{eqnarray*}
    P\left(\left.\tau_{\max}=0, I_{\max}=i, J(\tau_{\max})=j \right| (I(0),J(0))=(i_0,j_0)\right)
    \\
    & & \hspace{-1.5cm}= \ \delta_{(i,j),(i_0,j_0)} P_{(i_0,j_0)}(i_0),
\end{eqnarray*}
where $P_{(i_0,j_0)}(i)=P\left(\left.I_{\max}=i \right| (I(0),J(0))=(i_0,j_0)\right)$, and then derive iterative schemes for evaluating the \emph{restricted} Laplace-Stieltjes transform \begin{eqnarray*}
    \psi_{(i_0,j_0)}(\theta;i,j) &=& E\left[\left.e^{-\theta\tau_{\max}} 1\{I_{\max}=i,J(\tau_{\max})=j\} \right| (I(0),J(0))=(i_0,j_0)\right],
\end{eqnarray*}
for $\Re (\theta)\geq 0$ and states $(i,j)\in\cup_{k=i_0+1}^{\infty}l(k)$.
\par In terms of the \emph{taboo} Laplace-Stieltjes transform
\begin{eqnarray*}
    \varphi_{(i_0,j_0)}(\theta;i,j) &=& E\left[\left.e^{-\theta\tau_{l(i)}}1\{\tau_{l(i)}<\tau_{l(0)}, J(\tau_{l(i)})=j\} \right| (I(0),J(0))=(i_0,j_0)\right],
\end{eqnarray*}
for $\Re(\theta)\geq 0$, the restricted Laplace-Stieltjes transform of $\tau_{\max}$ on $\{\tau_{\max}>0,I_{\max}=i,J(\tau_{\max})=j\}$ is determined by
\begin{eqnarray}
    \psi_{(i_0,j_0)}(\theta;i,j) &=& \varphi_{(i_0,j_0)}(\theta;i,j)P_{(i,j)}(i),
    \label{eq:2}
\end{eqnarray}
for $\Re (\theta)\geq 0$, states $(i,j)\in\cup_{k=i_0+1}^{\infty}l(k)$ and the initial state $(i_0,j_0)\in {\cal S}\setminus l(0)$. To prove Eq.\ (\ref{eq:2}), we observe that $\psi_{(i_0,j_0)}(\theta;i,j)$ is linked to the time interval $[0,\tau_{\max})$ and the behavior of ${\cal X}$ during the residual interval $[\tau_{\max},\tau_{l(0)})$ during which the subset of states $\cup_{k=i+1}^{\infty} l(k)$ is avoided. More concretely, the random time $\tau_{\max}$ on $\{I_{\max}=i,J(\tau_{\max})=j\}$ is identically distributed to the first-passage time $\tau_{l(i)}$ on $\{J(\tau_{l(i)})=j\}$ under the taboo of states in $l(0)$, which yields the first term $\varphi_{(i_0,j_0)}(\theta;i,j)$ in the right-hand side of (\ref{eq:2}). The second term $P_{(i,j)}(i)$ in the right-hand side of (\ref{eq:2}) guarantees that, provided that $(I(\tau_{\max}),J(\tau_{\max}))=(i,j)$, the maximum level that ${\cal X}$ visits during the residual interval $[\tau_{\max},\tau_{l(0)})$ is not greater than the level variable at time $\tau_{\max}$. Then, the proof of (\ref{eq:2}) follows by applying the strong Markov property to the stopping time $\tau_{(i,j)}$ ---defined as the time of visit to state $(i, j)$--- and the natural filtration.
\subsection{Taboo Laplace-Stieltjes transforms of $\tau_{l(i)}$}
\label{Subect:2.2}
In deriving expressions for $\varphi_{(i_0,j_0)}(\theta;i,j)$, for states $(i,j)\in\cup_{k=i_0+1}^{\infty}l(k)$, we introduce taboo Laplace-Stieltjes transforms $\varphi_{(i',j')}(\theta;i,j)$ in a more general setting by replacing the initial state $(i_0,j_0)\in {\cal S}\setminus l(0)$ by the \emph{current} state $(i',j')$ of process ${\cal X}$ at any arbitrary time, for $(i',j')\in\cup_{k=1}^{i-1}l(k)$.
\\ \\
{\bf Theorem 1.}
{\it For any fixed state $(i,j)\in\cup_{k=i_0+1}^{\infty}l(k)$, the column vectors $\Phi_{i'}(\theta;i,j)=(\varphi_{(i',j')}(\theta;i,j) : j'\in\{0,...,M_{i'}\})$, for integers $i'\in\{1,...,i-1\}$, are uniquely characterized as the solution of the system of linear equations
\begin{eqnarray}
\left(\theta I_{1+M_{i'}}-Q_{i',i'}\right)\Phi_{i'}(\theta;i,j) &=& (1-\delta_{i',1})Q_{i',i'-1}\Phi_{i'-1}(\theta;i,j) \nonumber
\\
& & +Q_{i',i'+1}\Phi_{i'+1}(\theta;i,j), \quad i'\in\{1,...,i-2\},~~~ \label{eq:3}
\end{eqnarray}
with
\begin{eqnarray}
\left(\theta I_{1+M_{i-1}}-Q_{i-1,i-1}\right)\Phi_{i-1}(\theta;i,j) &=& Q_{i-1,i-2}\Phi_{i-2}(\theta;i,j)\nonumber
\\ & & +Q_{i-1,i}e_{1+M_i}(1+j). \label{eq:4}
\end{eqnarray}
Moreover, the column vectors $m_{i'}^{(n)}(i,j)$ consisting of the nth moments of $\tau_{l(i)}$ on $\{\tau_{l(i)}<\tau_{l(0)}, J(\tau_{l(i)})=j\}$, for $n\in\mathbb{N}$ and integers $i'\in\{1,...,i-1\}$, satisfy
\begin{eqnarray}
-Q_{i',i'}m_{i'}^{(n)}(i,j) &=& n m_{i'}^{(n-1)}(i,j) + (1-\delta_{i,1})Q_{i',i'-1}m_{i'-1}^{(n)}(i,j) \nonumber
\\
& & +Q_{i',i'+1}m_{i'+1}^{(n)}(i,j),\quad i'\in\{1,...,i-2\},
\label{eq:5}
\\
-Q_{i-1,i-1}m_{i-1}^{(n)}(i,j) &=& n m_{i-1}^{(n-1)}(i,j) + Q_{i-1,i-2}m_{i-2}^{(n)}(i,j),
\label{eq:6}
\end{eqnarray}
with $m_{i'}^{(0)}(i,j)=\left.\Phi_{i'}(\theta;i,j)\right|_{\theta=0}$.}
\\ \\
{\it Proof.}
  Let us fix the state $(i,j)\in\cup_{k=i_0+1}^{\infty}l(k)$. By conditioning on the first jump of process ${\cal X}$ when $(i',j')$ is its current state, it is found that
  \begin{description}
      \item[{\it (i)}] For states $(i',j')\in\cup_{k=1}^{i-2}l(k)$,
          \begin{eqnarray}
              \varphi_{(i',j')}(\theta;i,j) &=& (1-\delta_{i',1})\sum_{j''=0}^{M_{i'-1}}\frac{q_{(i',j'),(i'-1,j'')}}{\theta+q_{(i',j')}}\varphi_{(i'-1,j'')}(\theta;i,j) \nonumber
              \\
              & & +\sum_{j''=0, j''\neq j'}^{M_{i'}}\frac{q_{(i',j'),(i',j'')}}{\theta+q_{(i',j')}}\varphi_{(i',j'')}(\theta;i,j) \nonumber
              \\
              & & +\sum_{j''=0}^{M_{i'+1}}\frac{q_{(i',j'),(i'+1,j'')}}{\theta+q_{(i',j')}}\varphi_{(i'+1,j'')}(\theta;i,j),
              \label{eq:7}
          \end{eqnarray}
          where $q_{(i',j')}=-q_{(i',j'),(i',j')}$.
          \item[{\it (ii)}] For integers $j'\in\{0,...,M_{i-1}\}$,
          \begin{eqnarray}
              \varphi_{(i-1,j')}(\theta;i,j) &=& \sum_{j''=0}^{M_{i-2}}\frac{q_{(i-1,j'),(i-2,j'')}}{\theta+q_{(i-1,j')}}\varphi_{(i-2,j'')}(\theta;i,j) \nonumber
              \\
              & & +\sum_{j''=0, j''\neq j'}^{M_{i-1}}\frac{q_{(i-1,j'),(i-1,j'')}}{\theta+q_{(i-1,j')}}\varphi_{(i-1,j'')}(\theta;i,j) \nonumber
              \\
              & & +\frac{q_{(i-1,j'),(i,j)}}{\theta+q_{(i-1,j')}}.
              \label{eq:8}
          \end{eqnarray}
  \end{description}
  \par Straightforward algebra allows us to verify that, by multiplying (\ref{eq:7})-(\ref{eq:8}) by $\theta+q_{(i',j')}$ and $\theta+q_{(i-1,j')}$, respectively, the matrix versions of the resulting equations are (\ref{eq:3})-(\ref{eq:4}). Then, Eqs.\ (\ref{eq:5})-(\ref{eq:6}) are derived by taking derivatives in (\ref{eq:3})-(\ref{eq:4}) at $\theta=0$, and observing that the column vector $m_{i'}^{(n)}(i,j)=(-1)^n \left.d^n \Phi_{i'}(\theta;i,j)/d\theta^n\right|_{\theta=0}$ contains the taboo expectations
  \begin{eqnarray*}
      E\left[\left.\tau_{l(i)}^n 1\{\tau_{l(i)}<\tau_{l(0)}, J(\tau_{l(i)})=j\} \right| (I(0),J(0))=(i',j')\right],
  \end{eqnarray*}
  for any current state $(i',j')\in\cup_{k=1}^{i-1}l(k)$. This completes the proof. \ \ \ $\square$
\\
\par An important observation is that Eqs.\ (\ref{eq:3})-(\ref{eq:4}) (respectively, (\ref{eq:5})-(\ref{eq:6})) become a structured system of linear equations in the unknowns $\Phi_{i'}(\theta;i,j)$ (respectively, $m_{i'}^{(n)}(i,j)$) with $i'\in\{1,...,i-1\}$, whence its solution can be evaluated by applying any general-purpose procedure, such as block Gaussian elimination; see, e.g., Stewart \cite{Stewart1994}. This is summarized in Algorithms A-B, where the column vectors $m_{i'}^{(0)}(i,j)$ with
$i'\in\{1,...,i-1\}$ record the values $\left.\varphi_{(i',j')}(\theta;i,j)\right|_{\theta=0}$, for integers $j'\in\{0,...,M_{i'}\}$, and matrices $H_{i'}(i,j)$ with $i'\in\{1,...,i-1\}$ are defined by $\left.H_{i'}(\theta;i,j)\right|_{\theta=0}$; as a result, $m_{i'}^{(0)}(i,j)$ and $H_{i'}(i,j)$ in Algorithm B are derived from Algorithm A by selecting $\theta=0$.
\\ \\
{\small
  {\bf Algorithm A:} {\it Computation of the column vectors $\Phi_{i'}(\theta;i,j)$ with $i'\in\{1,...,i-1\}$, for $\Re(\theta)\geq 0$ and a fixed state $(i,j)\in\cup_{k=i_0+1}^{\infty}l(k)$.}
  \par Step 1: $i':=1$;
  \par \hspace{0.5cm} $H_{i'}(\theta;i,j):=(\theta I_{1+M_{i'}}-Q_{i',i'})^{-1}Q_{i',i'+1}$;
  \par \hspace{0.5cm} while $i'<i-2$, repeat
  \par \hspace{1.0cm} $i':=i'+1$;
  \par \hspace{1.0cm} $H_{i'}(\theta;i,j):=(\theta I_{1+M_{i'}}-Q_{i',i'}-Q_{i',i'-1}H_{i'-1}(\theta;i,j))^{-1}Q_{i',i'+1}$;
  \par \hspace{0.5cm} $H_{i'}(\theta;i,j):=H_{i'}(\theta;i,j)e_{1+M_{i'+1}}(1+j)$.
  \par Step 2: $\Phi_{i'}(\theta;i,j):=H_{i'}(\theta;i,j)$;
  \par \hspace{0.5cm} while $i'>1$, repeat
  \par \hspace{1.0cm} $i':=i'-1$;
  \par \hspace{1.0cm} $\Phi_{i'}(\theta;i,j):=H_{i'}(\theta;i,j)\Phi_{i'+1}(\theta;i,j)$.
\\ \\
 {\bf Algorithm B:} {\it Computation of the column vectors $m_{i'}^{(n)}(i,j)$ with $i'\in\{1,...,i-1\}$, for $n\in\mathbb{N}$ and a fixed state $(i,j)\in\cup_{k=i_0+1}^{\infty}l(k)$, in terms of $m_{i'}^{(n-1)}(i,j)$ with $i'\in\{1,...,i-1\}$.}
  \par Step 1: $i':=1$;
  \par \hspace{0.5cm} $h_{i'}(i,j):=-nQ_{i',i'}^{-1} m_{i'}^{(n-1)}(i,j)$;
  \par \hspace{0.5cm} while $i'<i-1$, repeat
  \par \hspace{1.0cm}  $i':=i'+1$;
  \par \hspace{1.0cm}  $h_{i'}(i,j):=(-Q_{i',i'}-Q_{i',i'-1}H_{i'-1}(i,j))^{-1}$;
  \par \hspace{1.0cm}  $h_{i'}(i,j):=h_{i'}(i,j)(Q_{i',i'-1}h_{i'-1}(i,j)+n m_{i'}^{(n-1)}(i,j))$.
  \par Step 2: $m_{i'}^{(n)}(i,j):=h_{i'}(i,j)$;
  \par \hspace{0.5cm} while $i'>1$, repeat
  \par \hspace{1.0cm} $i':=i'-1$;
  \par \hspace{1.0cm} $m_{i'}^{(n)}(i,j):=H_{i'}(i,j) m_{i'+1}^{(n)}(i,j)+h_{i'}(i,j)$.}
\\
\par From Algorithms A-B, the taboo Laplace-Stieltjes transform of $\tau_{l(i)}$ on $\{\tau_{l(i)}<\tau_{l(0)}, J(\tau_{l(i)})=j\}$ is easily obtained as
\begin{eqnarray*}
     \varphi_{(i_0,j_0)}(\theta;i,j) &=& e^T_{1+M_{i_0}}(1+j_0)\Phi_{i_0}(\theta;i,j),
\end{eqnarray*}
and the expectation $E\left[\left.(\tau_{\max})^n 1\{I_{\max}=i,J(\tau_{\max})=j\} \right| (I(0),J(0))=(i_0,j_0)\right]$ is evaluated as
\begin{eqnarray*}
  e^T_{1+M_{i_0}}(1+j_0)m_{i_0}^{(n)}(i,j)P_{(i,j)}(i),
\end{eqnarray*}
for any state $(i,j)\in\cup_{k=i_0+1}^{\infty}l(k)$.
\subsection{Marginal distribution of $I_{\max}$}
\label{Subect:2.3}
Using processes under a taboo (see, e.g., Section 2.1 in G\'{o}mez-Corral and L\'{o}pez Garc\'{\i}a \cite{GC2012}), the mass function $\{P_{(i_0,j_0)}(i) : i\in\mathbb{N}_0\}$ of $I_{\max}$ is found to be $P_{(i_0,j_0)}(i)=0$ if $i\in\{0,...,i_0-1\}$; otherwise, $P_{(i_0,j_0)}(i)=F_{\max}(i;i_0j_0)-(1-\delta_{i,i_0})F_{\max}(i-1;i_0,j_0)$, where the probability distribution function $F_{\max}(i;i_0,j_0)=P\left(\left.I_{\max}\leq i \right| (I(0),J(0))=(i_0,j_0)\right)$ has the form
\begin{eqnarray}
  F_{\max}(i;i_0,j_0) &=& e^T_{c(i)}((1-\delta_{i_0,1})c(i_0-1)+1+j_0)\left(-T^{-1}(i)t_0(i)\right),
  \label{eq:9}
\end{eqnarray}
for $i\in\{i_0,i_0+1,...\}$. In this expression, $c(i)=\sum_{i'=1}^i (1+M_{i'})$ is the cardinality of the subset of states $\cup_{k=1}^i l(k)$, the matrix $T(i)$ consists of infinitesimal rates $q_{(i,j),(i',j')}$, for states $(i,j),(i',j')\in\cup_{k=1}^i l(k)$, and $t_0(i)$ is the column vector
\begin{eqnarray*}
  t_0(i) &=& \left(\begin{array}{c}
      Q_{1,0}1_{1+M_0} \\ 0_{c(i)-(1+M_1)}
      \end{array}\right).
\end{eqnarray*}
It is worth noting that $-T^{-1}(i)$ in (\ref{eq:9}) is a matrix of expected total times spent in states of $\cup_{k=1}^i l(k)$, whence $-T^{-1}(i)$ is a non-negative matrix and $-T^{-1}(i)t_0(i)$ contains probabilities that the first visit of ${\cal X}$ to $l(0)$ occurs in a finite time, provided that process ${\cal X}$ does not visit states in the subset $\cup_{k=i+1}^{\infty} l(k)$, from any initial state in $\cup_{k=1}^i l(k)$.
\par Algorithm C derives the mass function of $I_{\max}$ in a recursive manner. Its proof mostly repeats arguments of Ref. \cite[Section 2.1]{GC2012}, whence we omit the details.
\\ \\
{\small
  {\bf Algorithm C:} {\it Computation of the non-null probabilities $P_{(i_0,j_0)}(i)$, for $i\in\{i_0,i_0+1,...\}$ and a predetermined initial state $(i_0,j_0)\in {\cal S}\setminus l(0)$.}
  \par Step 1: $i:=1$;
  \par \hspace{0.5cm} $-T^{-1}(i):=-Q_{i,i}^{-1}$;
  \par \hspace{0.5cm} if $i=i_0$, then
  \par \hspace{1.0cm} $F_{\max}(i;i_o,j_0):=e^T_{c(i)}(1+j_0)(-T^{-1}(i)t_0(i))$;
  \par \hspace{1.0cm} $P_{(i_0,j_0)}(i):=F_{\max}(i;i_0,j_0)$.
  \par Step 2: While $F_{\max}(i;i_0,j_0)<1-\varepsilon$, repeat
  \par \hspace{0.5cm} $i:=i+1$;
  \par \hspace{0.5cm} $A_{1,2}:=\left(\begin{array}{c} 0_{c(i-2)\times (1+M_i)} \\ Q_{i-1,i} \end{array}\right)$;
  \par \hspace{0.5cm} $A_{2,1}:=\left(0_{(1+M_i)\times c(i-2)}, Q_{i,i-1}\right)$;
  \par \hspace{0.5cm} $B_{2,2}:=\left(-Q_{i,i}-A_{2,1}(-T^{-1}(i-1))A_{1,2}\right)^{-1}$;
  \par \hspace{0.5cm} $B_{2,1}:=B_{2,2}A_{2,1}(-T^{-1}(i-1))$;
  \par \hspace{0.5cm} $B_{1,2}:=-T^{-1}(i-1)A_{1,2}B_{2,2}$;
  \par \hspace{0.5cm} $B_{1,1}:=-T^{-1}(i-1)\left(I_{c(i-1)}+A_{1,2}B_{2,1}\right)$;
  \par \hspace{0.5cm} $-T^{-1}(i):=\left(\begin{array}{cc} B_{1,1} & B_{1,2} \\ B_{2,1} & B_{2,2} \end{array}\right)$;
  \par \hspace{0.5cm} if $i\geq i_0$, then
  \par \hspace{1.0cm} $F_{\max}(i;i_0,j_0):=e^T_{c(i)}(c(i_0-1)+1+j_0)\left(-T^{-1}(i)t_0(i)\right)$;
  \par \hspace{1.0cm} $P_{(i_0,j_0)}(i):=F_{\max}(i;i_0,j_0)-F_{\max}(i-1;i_0,j_0)$.}
\\
\par Note that, from an appropriate use of Algorithms A and C, it is possible to evaluate the following expressions for the hitting probabilities at time $\tau_{\max}$ and the $n$th moment of $J(\tau_{\max})$ on $\{I_{\max}=i\}$:
\begin{description}
  \item[{\it (i)}] For states $(i,j)\in\cup_{k=i_0+1}^{\infty}l(k)$,
      \begin{eqnarray*}
          P\left(\left.I_{\max}=i,J(\tau_{\max})=j \right| (I(0),J(0))=(i_0,j_0)\right)
          \\
          & & \hspace{-1.5cm}= \ \left.\varphi_{(i_0,j_0)}(\theta;i,j)\right|_{\theta=0}P_{(i,j)}(i).
      \end{eqnarray*}
  \item[{\it (ii)}] For $n\in\mathbb{N}_0$ and integers $i\in\{i_0+1,i_0+2,...\}$,
      \begin{eqnarray*}
          E\left[\left.(J(\tau_{\max}))^n 1\{I_{\max}=i\} \right| (I(0),J(0))=(i_0,j_0)\right]
           \\
           & & \hspace{-2.5cm} = \ \sum_{j=0}^{M_i}j^n\left.\varphi_{(i_0,j_0)}(\theta;i,j)\right|_{\theta=0}P_{(i,j)}(i).
      \end{eqnarray*}
\end{description}
\par To conclude the section, we briefly present the solution for a finite LD-QBD process.
\\ \\
{\bf Remark 2.}
{\it For a process ${\cal X}$ on the finite space state ${\cal S}=\cup_{i=0}^Nl(i)$, Algorithms A-B should be executed progressively for any state $(i,j)\in\cup_{k=i_0+1}^{N}l(k)$ and a fixed initial state $(i_0,j_0)\in\cup_{k=1}^{N-1}l(k)$. In Algorithm C, the probabilities $P_{(i_0,j_0)}(i)$ can be evaluated for integers $i\in\{i_0,i_0+1,...,N\}$ by modifying Step 2 with the condition $i<N-1$ and noting that $F_{\max}(N;i_0,i_0)=1$.}
\section{Epidemic models}
\label{Sect:3}
In this section, the analytical solution in Section \ref{Sect:2} is exemplified with two epidemic models: the SIS model for a pathogen transmitted both vertically and horizontally; and the SIR model with constant population size. For convenience, the number $I(t)$ is related to infectious individuals in Subsections \ref{Subsect:3.1}-\ref{Subsect:3.2}, and the level of process ${\cal X}$ is suitably defined according to the extinction of the population (Subsection \ref{Subsect:3.1}) and the duration of an outbreak (Subsection \ref{Subsect:3.2}).
\subsection{The SIS epidemic model for vertically and horizontally transmitted pathogens}
\label{Subsect:3.1}
In the terminology of Reuter \cite{Reuter1961}, the SIS model with vertical and horizontal transmission is defined in Ref. \cite{GC2023b} as a competition process ${\cal Y}=\{(I(t),S(t)): t\geq 0\}$, where $I(t)$ and $S(t)$ are related to the number of infectious and susceptible individuals at time $t$, respectively. For states $(i,s), (i',s')\in\mathbb{N}_0\times\mathbb{N}_0$ with $(i',s')\neq (i,s)$, the dynamics of ${\cal Y}$ are specified by the infinitesimal rates
\begin{eqnarray}
  q_{(i,s),(i',s')} &=& \left\{\begin{array}{ll}
      \beta i s, & \hbox{if $i'=i+1$ and $s'=s-1$,}
      \\
      (1-p)\beta_I i, & \hbox{if $i'=i+1$ and $s'=s$,}
      \\
      \delta_S s, & \hbox{if $i'=i$ and $s'=s-1$,}
      \\
      p\beta_I i+\beta_S s, & \hbox{if $i'=i$ and $s'=s+1$,}
      \\
      \delta_I i, & \hbox{if $i'=i-1$ and $s'=s$,}
      \\
      \gamma i,  & \hbox{if $i'=i-1$ and $s'=s+1$,}
      \end{array}\right.
      \label{eq:SIR}
\end{eqnarray}
and $q_{(i,s),(i,s)}=-(\beta is+(\delta_I+\beta_I+\gamma)i+(\delta_S+\beta_S)s)$, for $(i,j)\in\mathbb{N}_0\times\mathbb{N}_0$, where $\beta$, $\beta_S$, $\beta_I$, $\delta_S$, $\delta_I$ and $\gamma$ are strictly positive, and $p\in (0,1)$. In particular, the value $\beta$ in (\ref{eq:SIR}) is the contact rate; parameters $\beta_S$ and $\delta_S$ (respectively, $\beta_I$ and $\delta_I$) are the birth and death rates of susceptible (respectively, infectious) individuals; $\gamma$ is the recovery rate of an infectious individual; and  $p$ represents the proportion of the offspring of infectious parents that is born as susceptible.
\par
\begin{figure}[h!]
\centering
\begin{center}
      \includegraphics[width=12cm]{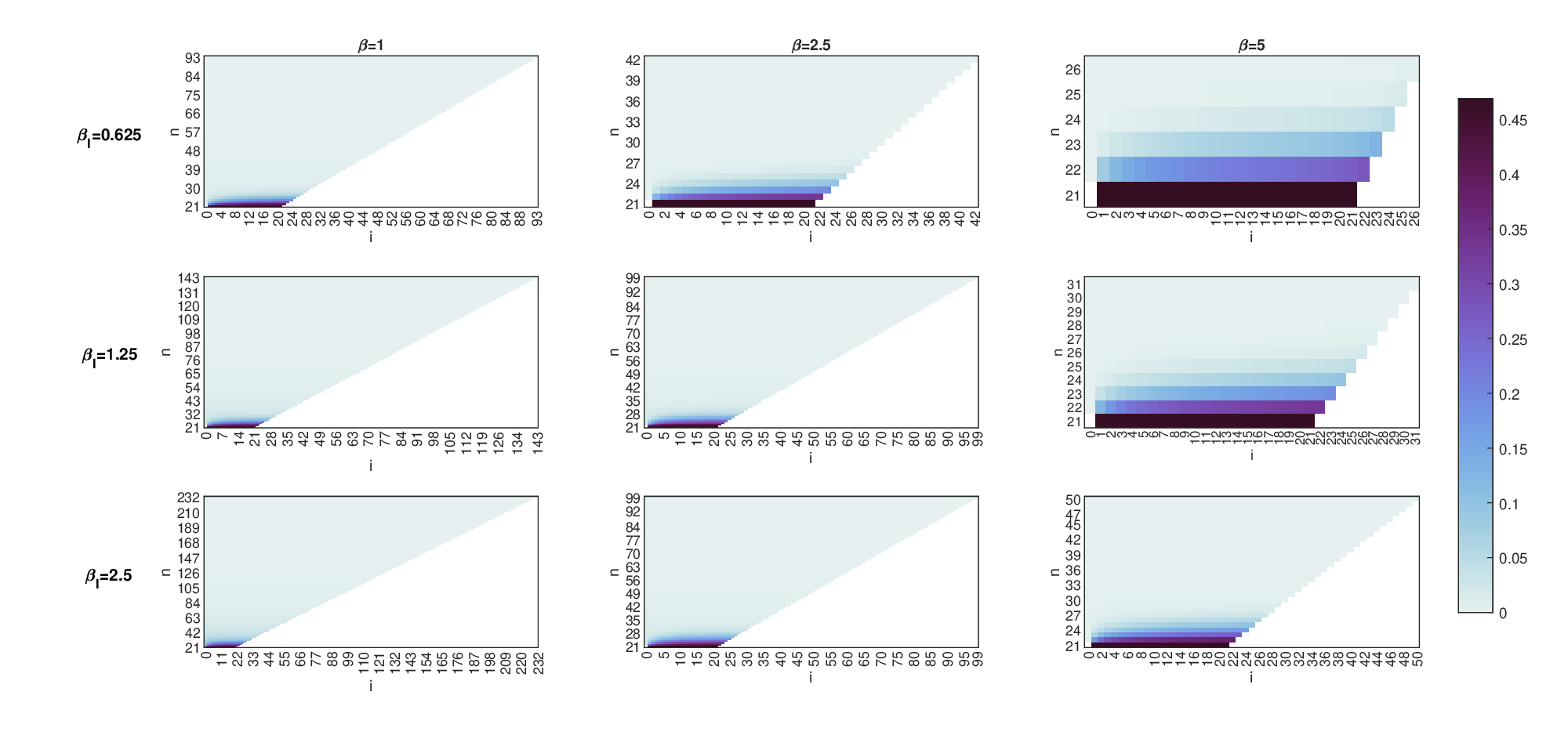}
\end{center}
\caption{Conditional probabilities $P\left(\left.N_{\max}=n, I(\tau_{\max})\leq i \right| (N(0),I(0))=(21,1)\right)$, for integers $n\in\{21,...,n_{0.95}\}$ and $i\in\{\delta_{n,21},...,n\}$, in scenarios with rates $\beta\in\{1.0,2.5,5.0\}$ and $\beta_I\in\{0.625,1.25,2.5\}$.}
\label{fig:SIS1}
\end{figure}

\begin{figure}
\centering
\begin{center}
      \includegraphics[width=12cm]{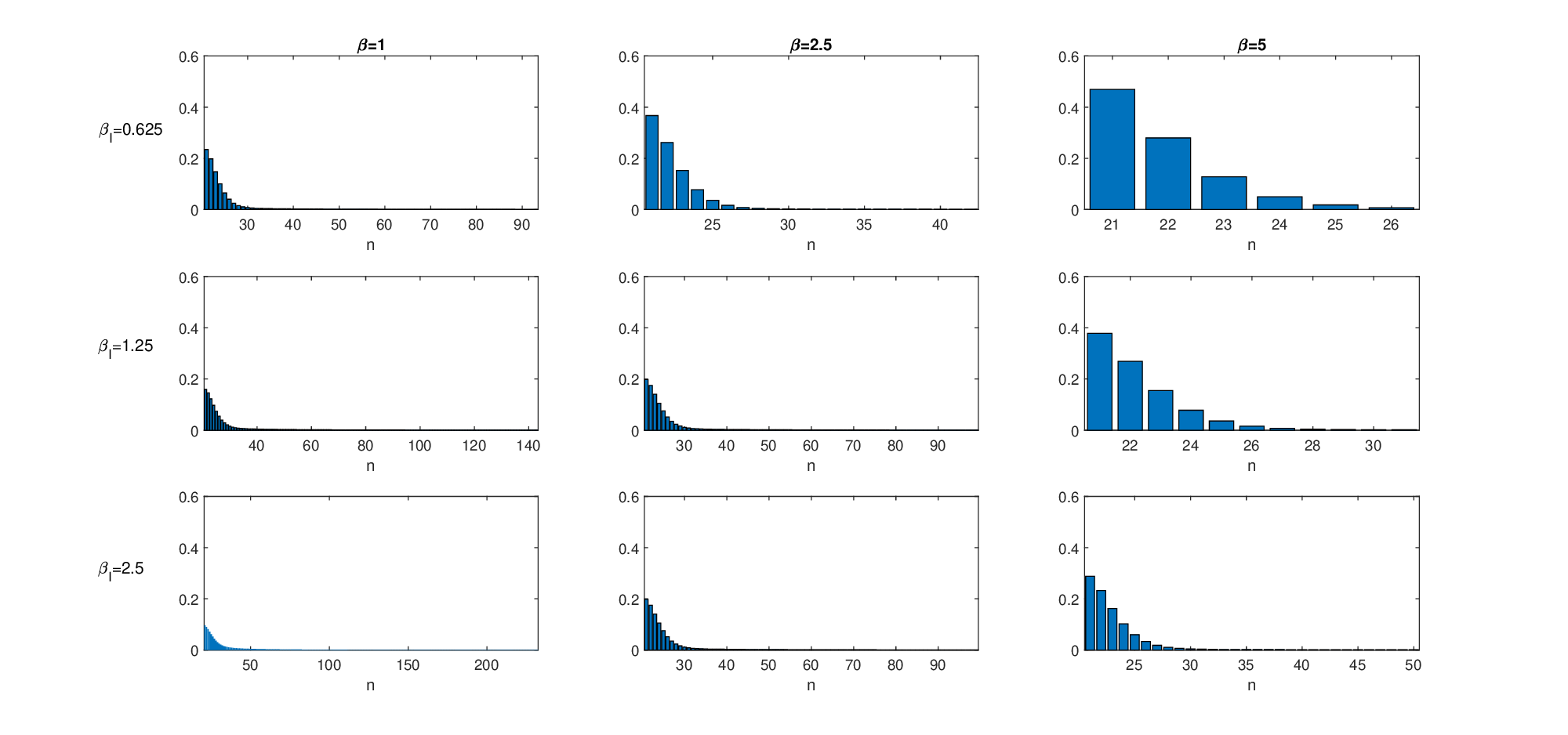}
\end{center}
\caption{The mass function $\{P\left(\left.N_{\max}=n \right| (N(0),I(0))=(21,1)\right): n\in\{21,...,n_{0.95}\}\}$ in scenarios with rates $\beta\in\{1.0,2.5,5.0\}$ and $\beta_I\in\{0.625,1.25,2.5\}$.}
\label{fig:SIS2}
\end{figure}

\begin{figure}
\centering
\begin{center}
      \includegraphics[width=12cm]{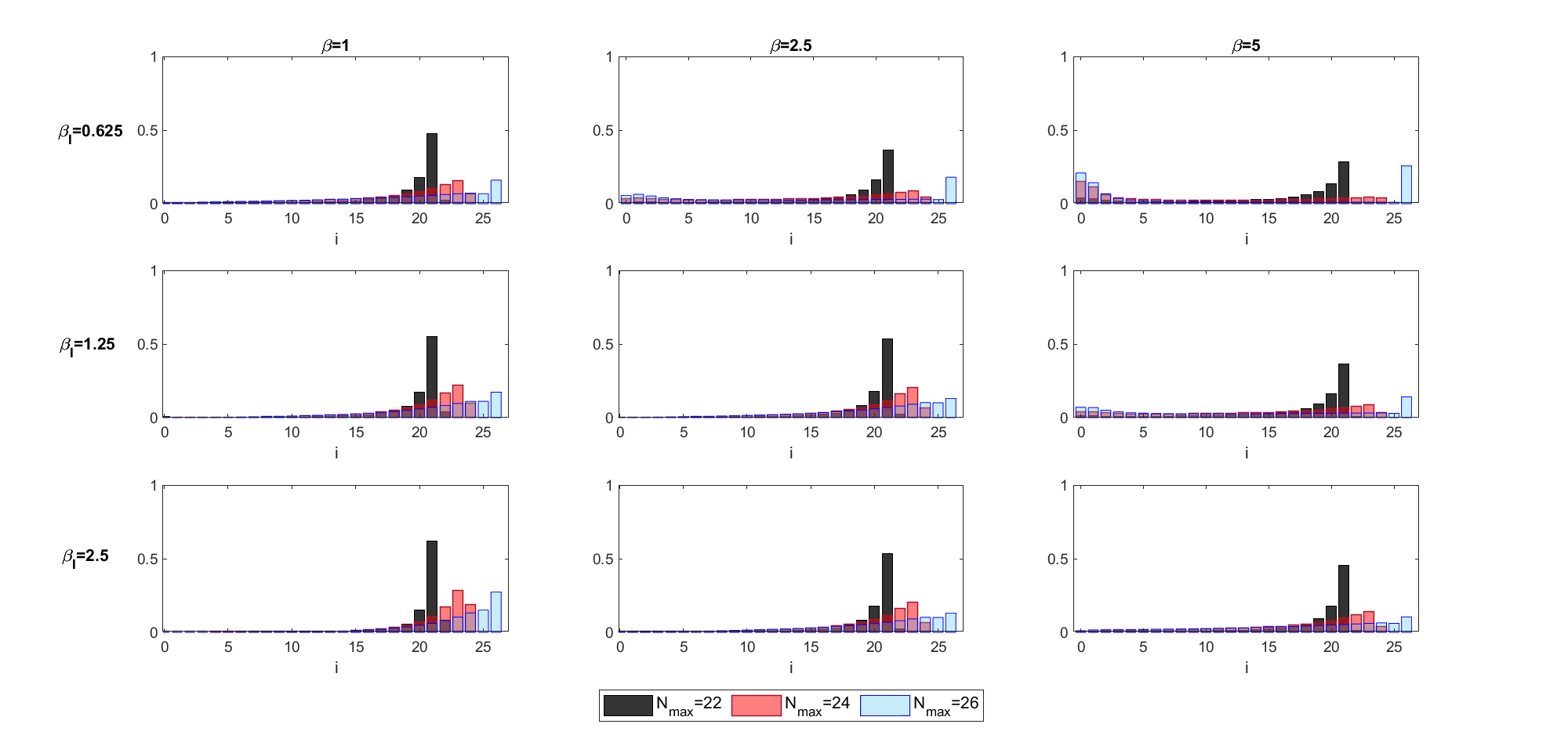}
\end{center}
\caption{Conditional probabilities
$\{P\left(\left.I(\tau_{\max})=i \right| (N(0),I(0))=(21,1), N_{\max}=n\right): i\in\{0,...,n\}\}$, for selected numbers $n\in\{22,24,26\}$, in scenarios with rates $\beta\in\{1.0,2.5,5.0\}$ and $\beta_I\in\{0.625,1.25,2.5\}$.}
\label{fig:SIS3}
\end{figure}

 In Figures \ref{fig:SIS1}-\ref{fig:SIS3}, we focus on the life cycle of a population initially consisting of a single infectious individual and twenty susceptible ones, i.e., $I(0)=1$ and $S(0)=20$. Thus, the analytical results in Section \ref{Sect:2} have to be applied here to process ${\cal X}=\{(N(t),I(t)): t\geq 0\}$ with $N(t)=I(t)+S(t)$, instead of ${\cal Y}$. Process ${\cal X}$ can be routinely formulated as a LD-QBD process\footnote{Specifically, we let $N(t)$ and $I(t)$ be the level and the phase variables of ${\cal X}$, respectively. Hence, $N_{\max}=N(\tau_{\max})$ is the maximum population size before extinction.} defined on ${\cal S}=\cup_{n=0}^{\infty}l(n)$ with levels $l(n)=\{(n,i): i\in\{0,...,n\}\}$, for $n\in\mathbb{N}_0$, in such a way that the time taken to access the absorbing state $(0,0)$ (i.e., level $l(0)$) corresponds to the extinction time. This allows quantifying the incidence of the disease at the time of maximum population size through the number $I(\tau_{\max})$ of infectious individuals, which acts as the phase of ${\cal X}$ at time $\tau_{\max}$.
\par Similarly to Ref. \cite[Section 4]{GC2023b}, birth and death rates in our experiments are assumed to be related to each other through the equalities $\delta_S=\beta_S$, $\beta_S=2\beta_I$ and $\delta_I=2\delta_S$. Hence, the compartment of susceptible individuals is expected to be demographically stable, offspring from infectious parents are expected to be born less frequently than those from susceptible parents, and susceptible individuals have a longer life expectancy than infectious ones. With the selections $p=0.2$ and $\gamma^{-1}=1.0$, it is supposed that $20\%$ of newborn offspring of infectious parents will be susceptible and that an infectious individual will be recovered, on average, in $1$ unit of time. Contact rates $\beta\in\{1.0,2.5,5.0\}$ and birth rates $\beta_I\in\{0.625,1.25,2.5\}$ yield nine scenarios in Figures \ref{fig:SIS1}-\ref{fig:SIS3}, for which the extinction of the population is certain and occurs in a finite mean time\footnote{Since process ${\cal X}$ is regular, the inequality $\min\{\beta_S,\beta_I\}<\delta_S+\delta_I$ can be seen as a sufficient condition for ${\cal X}$ to be absorbed almost surely, provided that $\min\{\delta_S,\delta_I\}=\min\{\beta_S,\beta_I\}$; see Ref. \cite[Theorems 1 and 3]{GC2023a}. Indeed, the absorption of ${\cal X}$ occurs in a finite mean time from any initial non-absorbing state.}. For convenience, intervals on the axes \emph{ox} and \emph{oy} in Figure \ref{fig:SIS1}, and on the \emph{ox} axis in Figures \ref{fig:SIS2}-\ref{fig:SIS3} are chosen to concentrate at least $0.95$ of the probability, whence the integer $n_{0.95}$ in these figures satisfies $F_{\max}(n_{0.95}-1;21,1)< 0.95 \leq F_{\max}(n_{0.95};21,1)$.
\par It is seen that, in our experiments, the transmission dynamics of the pathogen are limited and the duration of the outbreak is expected to be short. This is mainly related to the values of $P\left(\left.N_{\max}=n, I(\tau_{\max})\leq i \right| (N(0),I(0))=(21,1)\right)$ in Figure \ref{fig:SIS1}, especially because the most significant ones correspond to the integers $n=21$ and $i=1$ (i.e., the initial population size and the initially infectious individual), regardless of the choice of $(\beta,\beta_I)$. In Figure \ref{fig:SIS2}, it is seen that the mass function of $N_{\max}$ is unimodal and is mainly concentrated on the initial population size and the closest values to it. Note that, for a fixed value of $\beta_I$, shorter life cycles are expected to occur when isolation measures between compartments are relaxed (i.e., increasing values of $\beta$), which is explained by seeing that infectious individuals die at a higher rate than susceptible ones. Although a higher birth rate from infectious parents leads to an increased number of infectious newborns, it will also generate a greater number of susceptible newborns, who are born in a lower proportion than infectious ones but with longer life spans and a higher birth rate.
Therefore, heavier-tailed distributions of $N_{\max}$ ---which would predict longer durations of a life cycle--- are observed when isolation measures between compartments become more significant and, in particular, the birth rate from infectious parents increases.
\par A detailed description of the compartment of infectious individuals at the time of maximum population size is displayed in Figure \ref{fig:SIS3} in terms of the conditional probabilities $P\left(\left.I(\tau_{\max})=i \right| (N(0),I(0))=(21,1), N_{\max}=n\right)$, for integers $i\in\{0,...,n\}$ and selections $n\in\{22,24,26\}$. The unimodal conditional distribution of $I(\tau_{\max})$, provided that $N_{\max}=n$, reflects how a more accurate estimation of the maximum population size $N_{\max}$ will allow a better prediction of the incidence of the disease in terms of $I(\tau_{\max})$, as observed, for example, in the case $n=22$ and any of our scenarios. For some values of $(\beta,\beta_I)$, the selection $n=26$ leads to a bimodal conditional distribution of $I(\tau_{\max})$, with a first peak at a small integer $i$ and a second one at $i=n$. This behavior is linked to the inherent stochasticity of the Markov chain model, for which sample paths of process ${\cal X}$ on $\{N_{\max}=n\}$ in the case $n=26$ and some pairs $(\beta,\beta_I)$ may make it significantly plausible that the outbreak has a short duration in time and a minimal incidence in terms of the number $I(\tau_{\max})$ of infectious individuals, and at the same time a major outbreak with values of $I(\tau_{\max})$ comparable to the maximum population size can also be observed.
\subsection{The SIR epidemic model}
\label{Subsect:3.2}
In the SIR model (Kermack and McKendrick \cite{Kermack1927}; Neuts and Li \cite{Neuts1996}), the process ${\cal X}=\{(I(t),S(t)): t\geq 0\}$ records the number $I(t)$ of infectious individuals and the number $S(t)$ of susceptible individuals, whereas the number $R(t)$ of removed or recovered individuals is given by $R(t)=N-I(t)-S(t)$, where $N$ is the (constant) population size. The process ${\cal X}$ is then defined on the state space
\begin{eqnarray*}
  {\cal S} &=& \{(i,s) : i\in\{0,...,N_I+N_S-s\}, s\in\{0,...,N_S\}\},
\end{eqnarray*}
for initial numbers of $N_I$ infectious individuals and of $N_S$ susceptible individuals ---consequently, $N=N_I+N_S$ with $R(0)=0$---, and its $q$-matrix $Q$ has entries
\begin{eqnarray*}
  q_{(i,s),(i',s')} &=& \left\{\begin{array}{ll}
      N^{-1}\beta i s, & \hbox{if $i'=i+1$ and $s'=s-1$,}
      \\
      -(N^{-1}\beta is +\gamma i), & \hbox{if $i'=i$ and $s'=s$,}      \\
      \gamma i,  & \hbox{if $i'=i-1$ and $s'=s$,}
      \\
      0, & \hbox{otherwise,}
      \end{array}\right.
\end{eqnarray*}
for states $(i,s),(i',s')\in {\cal S}$, where the parameters $\beta$ and $\gamma$ are the contact and recovery rates, respectively, and are assumed to be strictly positive.
\begin{figure}
	\centering
	\begin{center}
        \includegraphics[width=10cm]{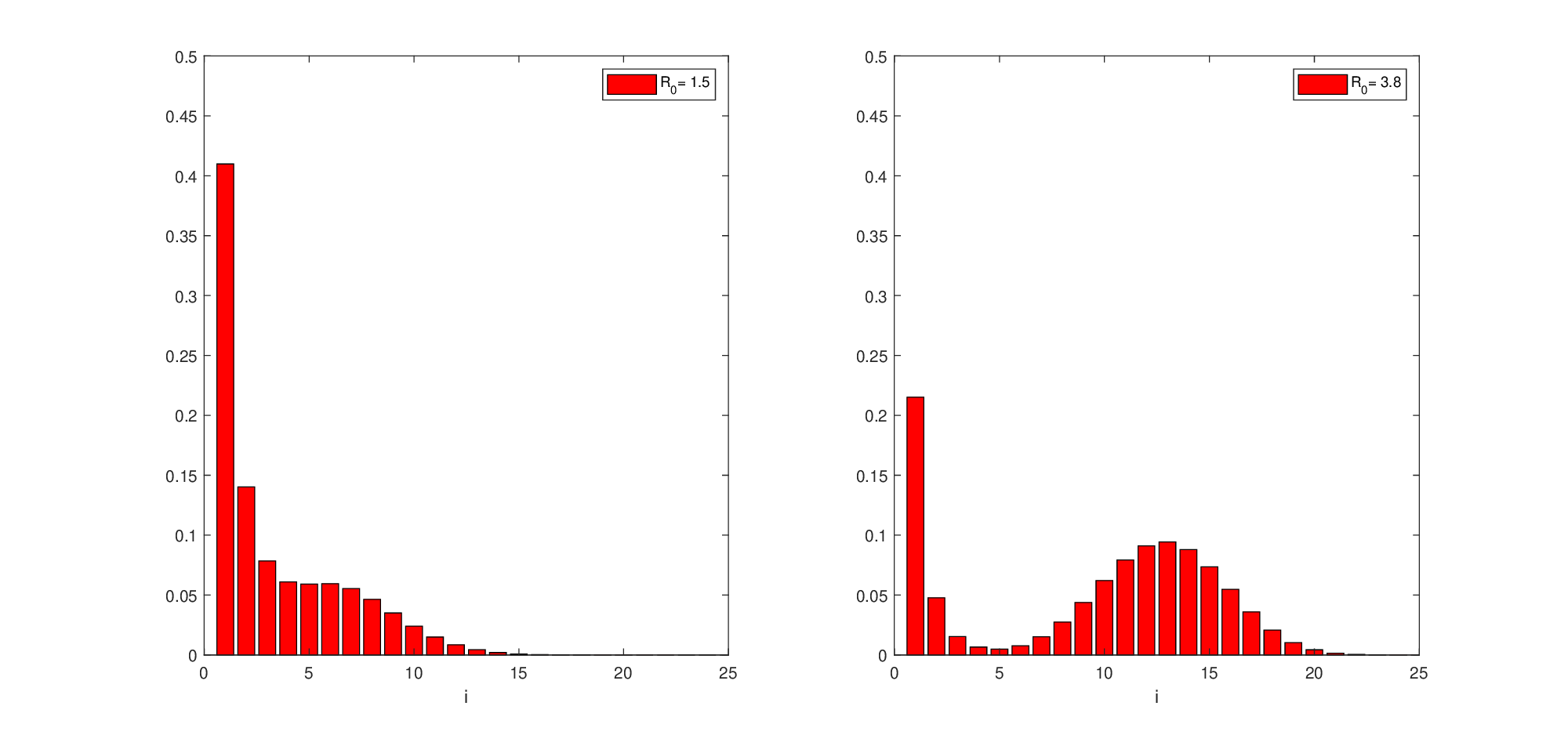}
	\end{center}
  \caption{ The mass function
  $\{P\left(\left.I_{\max}=i \right| (I(0),S(0))=(1,24)\right): i\in\{1,...,25\}\}$ for populations with value ${\cal R}_0=1.5$ (left) and $3.8$ (right). \hspace{2.5cm}}
	\label{fig:SIR1}
\end{figure}
\begin{figure}
  \centering
  \begin{center}
        \includegraphics[width=10cm]{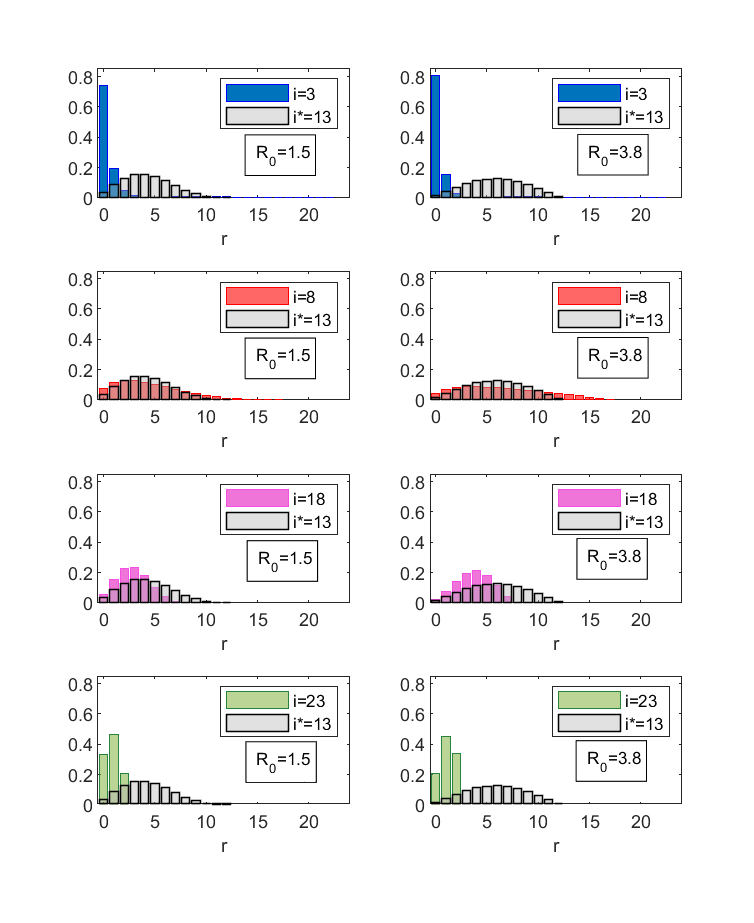}
  \end{center}
  \caption{Conditional probabilities
 $\{P\left(\left.R(\tau_{\max})=r \right| (I(0),S(0))=(1,24), I_{\max}=i\right): r\in\{0,...,25-i\}\}$, for selected numbers $i\in\{3,8,13,18,23\}$, and populations with value ${\cal R}_0=1.5$ (left) and $3.8$ (right).}
  \label{fig:SIR2}
\end{figure}
\par In Figures \ref{fig:SIR1}-\ref{fig:SIR2} and Table \ref{tab:SIR}, we consider a small community consisting initially of one infectious and twenty-four susceptible individuals, i.e., $I(0)=1$ and $S(0)=24$. An infectious individual is assumed to recover, on average, in $\gamma^{-1}=1$ unit of time, and the contact rate $\beta$ is determined from the basic reproductive number ${\cal R}_0=\gamma^{-1}\beta$. The choice of ${\cal R}_0$ values corresponds to the estimated mean values in Ref. \cite[Table 1]{VandenD2017} for the 1918 Spanish influenza. Figure \ref{fig:SIR1} illustrates the transmission dynamics in terms of the mass function\footnote{Since $I(t)$ and $S(t)$ specify the level and the phase variables of ${\cal X}$, respectively, $I_{\max}=I(\tau_{\max})$ represents the maximum size of the compartment of infectious individuals during an outbreak.} $\{P\left(\left.I_{\max}=i \right| (I(0),S(0))=(1,24)\right): i\in\{1,...,25\}\}$, which is plotted for values ${\cal R}_0=1.5$ (spring wave) and $3.8$ (fall wave), and Table \ref{tab:SIR} gives the expected values $E\left[\left.\tau_{\max} \right| (I(0),S(0))=(1,24), I_{\max}=i\right]$, for selected numbers $i\in\{3,8,13,18,23\}$. We notice that, not surprisingly, the less interaction between individuals, the less spread of the pathogen will be observed. In particular, the maximum number $I_{\max}$ of infectious individuals in the spring wave has an unimodal distribution with a clear peak at point $i=1$, meaning that the event $\{I_{\max}=1\}$ occurs more frequently than the others. As a result, the disease is more likely to be eradicated as soon as the initially infectious individual recovers than to spread to the subpopulation of susceptible individuals. However, the latter may also occur with significant probability. A more relevant interaction between individuals in the fall wave yields a bimodal distribution of $I_{\max}$, with a first peak at point $i=1$ and a second one at $i=13$. This illustrates how increasing values of ${\cal R}_0$ may cause major outbreaks ---i.e., outbreaks with a significant number of infectious cases--- at the same time that the disease could also disappear at an early stage of the epidemic. This fact highlights the need to analyze the maximum number $I_{\max}$ of infectious individuals as a random variable (i.e., in terms of its mass function) rather than its expected value or its deterministic counterpart.
\begin{table}
  \centering %
  \caption{Expected values $E\left[\left.\tau_{\max} \right| (I(0),S(0))=(1,24), I_{\max}=i\right]$ for selected numbers $i\in\{3,8,13,18,23\}$, and populations with ${\cal R}_0=1.5$ and $3.8$.}%
  \label{tab:SIR}
  \begin{tabular}{rrc}
            \\ \hline
            ${\cal R}_0$ & $i$ & $E\left[\left.\tau_{\max} \right| (I(0),S(0))=(1,24), I_{\max}=i\right]$ \\
            \hline
            1.5 &  3 & 0.84726 \\
                &  8 & 2.44037 \\
                & 13 & 2.27654 \\
                & 18 & 2.07524 \\
                & 23 & 1.91487 \\
            3.8 &  3 & 0.41006 \\
                &  8 & 1.73395 \\
                & 13 & 1.45888 \\
                & 18 & 1.26555 \\
                & 23 & 1.13595 \\
            \hline
    \end{tabular}
  \end{table}
\par To further investigate how the values of $I_{\max}$ affect the subpopulation of recovered individuals at the time $\tau_{\max}$ of greatest incidence of the pathogen, a detailed description of $R(\tau_{\max})$ is displayed in Figure \ref{fig:SIR2}, in terms of the mass function $\{P\left(\left.R(\tau_{\max})=r \right| (I(0),S(0))=(1,24), I_{\max}=i\right): r\in\{0,...,25-i\}\}$, for integers $i\in\{3,8,13,18,23\}$. It is important to note that the most likely event $\{I_{\max}=i\}$ in Figure \ref{fig:SIR1} is linked to the choice $i=1$, irrespectively of ${\cal R}_0$, which leads to $R(\tau_{\max})=0$ almost surely (since $\tau_{\max}=0$), provided that $(I(0),S(0))=(1,24)$ and $I_{\max}=1$. The particular selection $i^*=13$ is related to the second peak of the bimodal distribution of $I_{\max}$ in the fall wave. Therefore, the conditional probabilities $\{P\left(\left.R(\tau_{\max})=r \right| (I(0),S(0))=(1,24), I_{\max}=i^*\right): r\in\{0,...,25-i^*\}\}$ must be seen as an accurate description of the variability of $R(\tau_{\max})$ for a major outbreak in the case ${\cal R}_0=3.8$, whereas they should not be usually appropriate when ${\cal R}_0=1.5$ as $i^*=13$ leads to an event $\{I_{\max}=i^*\}$ that it will hardly ever occur. Our experiments always found that the conditional distribution of $R(\tau_{\max})$ is unimodal, regardless of the selected integer $i$.
\section{Conclusions}
\label{Sect:4}
For a regular LD-QBD process ${\cal X}=\{(I(t),J(t)): t\geq 0\}$ with countably many states, the probability law of the random vector $(\tau_{\max},I_{\max},J(\tau_{\max}))$ has been characterized under the assumption that the first passage from any initial state $(i_0,j_0)\in {\cal S}\setminus l(0)$ to level $l(0)$ is certain. For this purpose, efficient algorithms for computing the marginal distribution of $I_{\max}$, and restricted Laplace-Stieltjes transforms of $\tau_{\max}$ on the sample paths of ${\cal X}$ such that $\{I_{\max}=i,J(\tau_{\max})=j\}$, have been developed as an extension of the well-known block-Gaussian elimination technique, whence the algorithmic complexity of the solution in Subsections \ref{Subect:2.2} and \ref{Subect:2.3} is similar to that in Refs. \cite[Algorithm A]{Gaver1984} and \cite[Algorithms 1A-3A]{GC2018} in the finite case.
\par The computational approach can be thought of as an alternative one to that of Mandjes and Taylor \cite{Mandjes2016}, and Javier and Fralix \cite{Javier2023} for the running maximum $\bar{I}(T)$ attained by process ${\cal X}$ during a time interval $[0,T]$ when $T$ is either independent and exponentially distributed time, or a fixed time. Here, time $T$ is described as the first time that process ${\cal X}$ visits states in level $l(0)$, so $T$ is closely linked to the dynamics of ${\cal X}$ and, for practical use, amounts to the duration of a life cycle in a population, of an outbreak in epidemics, and of a busy period in a queueing model. In the context of epidemics, an aspect that deserves further exploration is how phase $J(\tau_{\max})$ and its probability law ---when level $I(t)$ in process ${\cal X}$ is defined as the number of infectious individuals--- could be used to characterize a stochastic version of the herd immunity threshold (see, e.g., Britton et al. \cite{Britton2020}), for which it is crucial to notice that no more endemic transmissions are expected to occur after time $\tau_{\max}$.
\section*{Acknowledgements}
This research was supported by Ministerio de Ciencia e Innovaci\'{o}n (Project PID2021-125871NB-B-I00), as well as by GNCS-INdAM and ``European Union -- Next Generation EU'' PRIN 2022 PNRR (Project P2022XSF5H). D. Taipe also acknowledges the support of Banco Santander-Universidad Complutense de Madrid (Pre-doctoral Researcher Contract CT63/19-CT64/19).
\section*{Data availability}
No data was used for the research described in the article. Details on numerical experiments in Section \ref{Sect:3} will be made available from the corresponding author upon reasonable request.
\section*{Authors' contributions}
All the authors equally contributed to this work.

\end{document}